\documentclass[11pt]{article}
\usepackage[T1]{fontenc}
\usepackage{lmodern,amsmath,amssymb,amsthm}
\usepackage[margin=1in]{geometry}
\usepackage[hidelinks]{hyperref}
\newtheorem{theorem}{Theorem}[section]
\newtheorem{lemma}[theorem]{Lemma}
\newtheorem{corollary}[theorem]{Corollary}
\newcommand{\Z}{\mathbb Z}
\newcommand{\F}{\mathbb F}
\newcommand{\End}{\operatorname{End}}
\newcommand{\tr}{\operatorname{tr}}
\newcommand{\ad}{\operatorname{ad}}
\newcommand{\SL}{\operatorname{SL}}

\newcommand{\slt}{\mathfrak{sl}_2}
\title{A congruence-group counterexample to Kourovka Problem 21.7}
\author{L\"u Gong\\[3pt]
  \small School of Mathematics and Statistics, Nantong University\\
  \small Jiangsu 226019, China\\
  \small\href{mailto:lieningzai1917@126.com}{\nolinkurl{lieningzai1917@126.com}}}
\date{}
\begin{document}
\maketitle
\begin{abstract}
For every prime $p>24$, let $H$ be the kernel of
$\SL_2(\Z/p^9\Z)\to\SL_2(\F_p)$. We prove that $H$ is not IYB,
whereas $H\times C_{p^8}$ is IYB. The proof excludes a pre-Lie
multiplication on the rank-three Lie ring with bracket $p[x,y]$ over
$\Z/p^8\Z$ using three trace identities. It follows that the supersoluble
group $H\times(C_{p^8}\rtimes C_2)$, with the involution acting by
inversion, has IYB Sylow subgroups and is not IYB. Its nilpotent residual
is cyclic and its Sylow $p$-subgroup has nilpotency class eight.
For $p=31$, an action of order three gives a supersoluble
counterexample of odd order $3\cdot31^{32}$. The same construction
gives an IYB $p$-group with a specified central cyclic direct factor
such that the quotient of the whole group by the unique subgroup of order $p$
in that factor is non-IYB.
\end{abstract}

\section{The finite groups and the question}
A finite group $P$ is an involutive Yang--Baxter (IYB) group if there is an abelian $P$-module $A$ and a
bijection $b:P\to A$ satisfying $b(xy)=b(x)+x\cdot b(y)$. Equivalently,
$P$ acts regularly by affine transformations on $A$, or $P$ is the
multiplicative group of a left brace with abelian additive group.
Every finite IYB group is soluble, and its Sylow subgroups are IYB
\cite{CJR}. Vendramin recorded a question of Rump asking whether a
finite non-IYB group could have all its Sylow subgroups IYB
\cite[Problem~2.35]{Ven}. That formulation allows $A_5$, whose Sylow
subgroups are $C_2\times C_2$, $C_3$ and $C_5$ up to isomorphism.
Ced\'o and Okni\'nski asked whether a finite soluble group must be IYB
whenever all its Sylow subgroups are IYB \cite[Question~1.1]{CO}. This is also the
formulation recorded by Ballester-Bolinches as Problem~21.7 of the
Kourovka Notebook \cite{KN}.

Bachiller constructed a finite $p$-group that admits no IYB structure
\cite{Bac}. Such a group already fails the Sylow hypothesis, since it
is its own Sylow $p$-subgroup. We construct a non-IYB $p$-group $H$
for which adjoining a cyclic direct factor gives an IYB group. An
involution acting by inversion on that factor has fixed subgroup $H$.
Denote by
$\operatorname{Dih}(C_m)$ the group $C_m\rtimes C_2$ defined by
inversion, of order $2m$.

\begin{theorem}\label{thm:main}
For a prime $p>24$, put $O=\Z/p^9\Z$ and
\[
 H=\ker\bigl(\SL_2(O)\longrightarrow\SL_2(\F_p)\bigr),\qquad
 G=H\times\operatorname{Dih}(C_{p^8}).
\]
Then $|H|=p^{24}$ and $|G|=2p^{32}$. The group $G$ is supersoluble,
and its nilpotent residual is the rotation subgroup $C_{p^8}$.
The Sylow $p$-subgroup $H\times C_{p^8}$ has nilpotency class eight
and is IYB. All Sylow $2$-subgroups of $G$ are IYB.
Both $H$ and $G$ are non-IYB.
\end{theorem}
Meng, Ballester-Bolinches and Jim\'enez-Seral proved that a finite
soluble group with Sylow subgroups of nilpotency class at most two is
IYB if its nilpotent residual has no section isomorphic to the
quaternion group $Q_8$ \cite[Theorem~A]{MBJ}. The cyclic residual in
Theorem~\ref{thm:main} satisfies this condition. The class-eight
Sylow subgroup shows that the class bound in their theorem cannot
be replaced by the assumption that each Sylow subgroup is IYB.
For groups of odd order, Ced\'o and Okni\'nski proved that the same
class-two bound on the Sylow subgroups suffices \cite[Theorem~2.5]{CO}.
Corollary~\ref{cor:odd} gives an odd-order supersoluble counterexample
when one assumes only the IYB property of the Sylow subgroups.

\section{A trace obstruction over a finite local ring}
Let $p>3$, $k\ge8$, $R=\Z/p^k\Z$, and let $V$ have basis $h,e,f$
over $R$. Write $[\cdot,\cdot]_0$ for the bracket specified by
$[h,e]_0=2e$, $[h,f]_0=-2f$, and $[e,f]_0=h$.
For a pre-Lie product $x*y$, the operators $L_x(y)=x*y$ satisfy
$[L_x,L_y]=L_{x*y-y*x}$. A pre-Lie structure on the Lie ring with
bracket $p[\cdot,\cdot]_0$ would therefore give the following two identities.
Every biadditive product here is $R$-bilinear.

\begin{lemma}\label{lem:obstruction}
There is no $R$-linear assignment $x\mapsto L_x\in\End_R(V)$ satisfying
\begin{equation}\label{eq:two}
 [L_x,L_y]=pL_{[x,y]_0},\qquad L_x(y)-L_y(x)=p[x,y]_0
 \qquad(x,y\in V).
\end{equation}
\end{lemma}
\begin{proof}
Write $D=L_h$, $E=L_e$, $F=L_f$ and $R'=\Z/p^{k-1}\Z$.
All identities with a factor $p$ below are first taken over $R$.
Cancelling that factor gives an identity over $R'$. Since $h,e,f$
are unit multiples of brackets, taking traces in the first part of
\eqref{eq:two} gives $\tr L_x=0$ over $R'$ for all $x$.

The tensor
\[
 B=h\otimes h+2e\otimes f+2f\otimes e=\sum_i a_i\otimes b_i
\]
is invariant under the adjoint action for the unscaled bracket. Indeed
$h$ gives zero on each of the three terms by weights. For $e$, the
contribution from $h\otimes h$ is $-2(e\otimes h+h\otimes e)$
and the other terms cancel it. For $f$ the two contributions are
$2(f\otimes h+h\otimes f)$ and its negative.
Set
\[
 w=\sum_iL_{a_i}b_i,qquad C=\sum_iL_{a_i}L_{b_i}=D^2+2EF+2FE.
\]
Using both identities in \eqref{eq:two} and then the invariance of $B$
gives, over $R$,
\[
 L_xw=Cx+p\sum_i\bigl(L_{[x,a_i]_0}b_i+L_{a_i}[x,b_i]_0\bigr)=Cx.
\]
Thus $C=L_w-p\ad_0(w)$, where $\ad_0(w)(x)=[w,x]_0$. The unscaled adjoint operators are traceless,
so $\tr C=0$ over $R'$. Cyclicity of the trace over $R$ gives
\[
 p\tr(D^2)=\tr([E,F]D)=\tr(E[F,D])=2p\tr(EF).
\]
Over $R'$ it follows that $0=\tr C=3\tr(D^2)$, so $\tr(D^2)=0$.

For the cubic trace, over $R$ again,
\[
 \begin{aligned}
 p\tr(D^3)&=\tr([E,F]D^2)=\tr(E[F,D^2])
             =2p\tr\bigl(E(FD+DF)\bigr),\\
 0&=\tr([E,F^2]E)=p\tr\bigl(E(DF+FD)\bigr).
 \end{aligned}
\]
The second line uses $[E,F^2]=p(DF+FD)$ and cyclicity of the trace.
Hence $\tr(D^3)=0$ over $R'$.

The Newton identities, with $2$ and $3$ invertible, give zero for all
three nonleading coefficients of the characteristic polynomial of $D$
over $R'$. Cayley--Hamilton now gives $D^3=0$ there. Consequently
$(\ad D)^5=0$, since each term in $(\ad D)^5(X)$ has at least three
copies of $D$ on one side of $X$. Applying this to $E$ and $F$ and using
$[D,E]=2pE$ and $[D,F]=-2pF$ yields $p^5E=p^5F=0$ over $R'$.
The second part of \eqref{eq:two} says $Ef-Fe=ph$, so $p^6h=0$
over $R'$. This contradicts $k-1\ge7$ and the fact that $h$ is a
basis vector.
\end{proof}

\section{Passage to finite groups}
We use the necessary direction of the finite-brace Lazard correspondence
\cite[Corollary~8.7.5]{TrThesis}, see also \cite{Tr}. We recall a
proof of this direction in order to keep track of the order bound.
For a finite $p$-group of class less than $p$, write $\log P$ for
its Lie ring under the classical Lazard correspondence.

\begin{lemma}\label{lem:lazard}
If $P$ is an IYB group of order $p^n$ with $n<p$, then $\log P$ admits a
biadditive pre-Lie product inducing its Lie bracket.
\end{lemma}
\begin{proof}
Take a regular affine representation of $P$ on an abelian group $A$.
Its group of linear parts preserves a series
$0=A_0<A_1<\cdots<A_n=A$ with all quotients of order $p$ on which it
acts trivially. To obtain the series, take a fixed line in the
$p$-torsion of each successive quotient using the $p$-group action.
Let $p^a$ be the exponent of $A$. On $A\oplus\Z/p^a\Z$ take the
associative ring $J$ of block endomorphisms
$\left(\begin{smallmatrix}T&v\\0&0\end{smallmatrix}\right)$ with
$T(A_i)\subseteq A_{i-1}$ and $v\in A$. We have $J^{n+1}=0$, and
the affine representation lands in $1+J$. Since $n<p$, logarithm
and exponential identify the Lazard Lie ring of $1+J$ with $J$ with
the commutator bracket. Functoriality gives a Lie homomorphism
$x\mapsto (T_x,d(x))$ from $\log P$ to this affine Lie ring.
The additive map $d$ is injective. If $d(x)=0$, the homogeneous
matrix for the exponential fixes $(0,1)$, and regularity forces
$x=0$. It is an additive isomorphism by cardinality. The bracket
identity gives
\[
 [T_x,T_y]=T_{[x,y]},\qquad d([x,y])=T_xd(y)-T_yd(x).
\]
Taking $L_x=d^{-1}T_xd$ gives the required product.
\end{proof}

\begin{proof}[Proof of Theorem~\ref{thm:main}]
Let $J=pM_2(O)$. Its ninth power is zero and $8<p$. The exponential
and logarithm series through degree eight give inverse maps
$J\leftrightarrow 1+J$. The logarithm identifies the group law with
the Baker--Campbell--Hausdorff law. In these coordinates, $H$ consists precisely of the
trace-zero matrices in $J$. For instance, lift a matrix in $J$ to
$pM_2(\Z_p)$ and reduce the identity
$\det(\exp X)=\exp(\tr X)$ modulo $p^9$. The scalar exponential on
$pO$ is injective. Multiplication by $p$ identifies
$\slt(\Z/p^8\Z)$ additively with the trace-zero matrices in $J$ and
turns their Lie bracket into $p[\cdot,\cdot]_0$ on $\slt(\Z/p^8\Z)$.
In particular $|H|=p^{24}$. To compute its nilpotency class, write
$V=\slt(\Z/p^8\Z)$ and $L=(V,p[\cdot,\cdot]_0)$. The basis brackets
give $[V,V]_0=V$, so induction in the lower central series gives
$\gamma_i(L)=p^{i-1}V$ for $i\ge1$. Thus $\gamma_8(L)=p^7V\ne0$
and $\gamma_9(L)=0$. The Lazard correspondence preserves nilpotency
class, so $H$ has class eight.

If $H$ were IYB, Lemma~\ref{lem:lazard} would apply since $24<p$.
It would supply the operators excluded by Lemma~\ref{lem:obstruction}
with $k=8$. Hence $H$ is non-IYB.

The group $P=1+J$ is IYB. Its additive module is $J$, the action is
left matrix multiplication by $g\in P$, and $b(g)=g-1$ is a
bijective cocycle. Every $g\in P$ can be written uniquely as $s h$
with $s$ scalar in $1+pO$ and $h\in H$. In fact $s$ is determined by
$s^2=\det g$, and squaring is an automorphism of the abelian
$p$-group $1+pO$. This scalar group is cyclic of order $p^8$, as
is seen from the scalar logarithm. Consequently $P\cong H\times C_{p^8}$,
and $P$ has the same nilpotency class as $H$.

Let $\alpha$ act trivially on $H$ and by inversion on $C_{p^8}$.
The fixed subgroup of this automorphism of $P$ is $H$.
Now $G=P\rtimes\langle\alpha\rangle$ has normal Sylow $p$-subgroup
$P$ and Sylow $2$-subgroups of order two. It is supersoluble, since
it is a direct product of a $p$-group and a dihedral group with cyclic
rotation group of prime-power order. Put $M=C_{p^8}$. The quotient
$G/M\cong H\times C_2$ is nilpotent. In any nilpotent quotient of
$G$, the images of $M$ and $\langle\alpha\rangle$ commute, since
they have coprime orders. Every such quotient therefore kills
$[M,\alpha]=M$, where the equality follows from inversion and the
fact that $p$ is odd. Hence $M$ is the nilpotent residual of $G$.

If $G$ were IYB, its normal Hall
subgroup $P$ would admit an $\alpha$-equivariant bijective cocycle
by \cite[Proposition~2.3]{Eis}. Such a cocycle maps $C_P(\alpha)$
bijectively onto the additive fixed subgroup. The latter is a
$C_P(\alpha)$-module by compatibility, so restriction would give an
IYB structure on $H$. This is impossible.
\end{proof}

\begin{corollary}\label{cor:odd}
Let $H$ be the group of Theorem~\ref{thm:main} for $p=31$, and write
$M=\Z/31^8\Z$ additively. There is an action of $C_3$ on $M$ such
that
\[
 G_{\mathrm{odd}}=H\times(M\rtimes C_3)
\]
is a supersoluble group of order $3\cdot31^{32}$ with IYB Sylow
subgroups. The group $G_{\mathrm{odd}}$ is non-IYB.
\end{corollary}
\begin{proof}
The polynomial $X^2+X+1$ has a simple root $5$ modulo $31$, since
$5^2+5+1=31$ and $2\cdot5+1=11$. This root lifts to a root $u$
modulo $31^8$ by Hensel's lemma. Then $u^3=1$ and $u-1$ is a unit,
so multiplication by $u$ has order three and fixes only zero in $M$.
Let a generator $\beta$ of $C_3$ act this way on $M$ and trivially
on $H$. The Sylow $31$-subgroup $P=H\times M$ is IYB by
Theorem~\ref{thm:main}, and the Sylow $3$-subgroups are cyclic.
Since $|H|=31^{24}$, the asserted order follows.
The subgroups $31^iM$ are normal in $G_{\mathrm{odd}}$, and successive
terms have index $31$. This series extends to a normal series with
prime-order factors because $G_{\mathrm{odd}}/M\cong H\times C_3$ is
nilpotent. Thus $G_{\mathrm{odd}}$ is supersoluble.

If $G_{\mathrm{odd}}$ were IYB, its normal Hall subgroup $P$ would
admit a $\langle\beta\rangle$-equivariant bijective cocycle
\cite[Proposition~2.3]{Eis}. Restriction to $C_P(\beta)=H$ would
give a bijective cocycle onto the fixed subgroup of the additive
module, contradicting Theorem~\ref{thm:main}.
\end{proof}

\section{Removing a single central subgroup}
\begin{corollary}\label{cor:peel}
For every prime $p>24$, there is an IYB group $S=M\times H$ with a
specified nontrivial central cyclic $p$-group $M$ such that $S/U$ is
non-IYB for the unique subgroup $U$ of order $p$ in $M$. One can
take the group $H$ from Theorem~\ref{thm:main} and $|M|=p^j$ for some
$1\le j\le6$.
\end{corollary}
\begin{proof}
The groups $S_i=H\times C_{p^i}$ start with $S_0=H$ non-IYB and
include the IYB group $S_8=P$. The least index $j$ for which $S_j$
is IYB is positive. Its specified cyclic factor has only one subgroup
of order $p$, and quotienting by this subgroup gives $S_{j-1}$.

For the stated upper bound, take the subring
$J_0=\{X\in J:\tr X\in p^2O\}$.
It contains all products of two elements of $J$.
Since $\det(1+X)=1+\tr X+\det X$ and $\det X\in p^2O$,
$1+J_0$ is the subgroup of $P$ with determinant in $1+p^2O$.
Taking its scalar square roots gives $1+J_0\cong H\times C_{p^7}$.
The additive subgroup of scalar matrices $p^8O I_2$ is an ideal annihilating
$J_0$ on both sides. The adjoint group of the quotient ring
$J_0/(p^8O I_2)$ is therefore IYB and is isomorphic to
$H\times C_{p^6}$. This proves $j\le6$.
\end{proof}
This argument does not determine $j$. In particular, the question
whether every finite $p$-group $K$ for which $K\times C_p$ is IYB must itself be IYB remains open in this argument. The explicit group $G$ in
Theorem~\ref{thm:main} is independent of $j$.

\begin{corollary}
Let $p>24$ be prime. For every prime $q\ne p$ there is an IYB
$p$-group admitting an automorphism of order $q$ that preserves no
addition making the group a left brace with abelian additive group.
\end{corollary}
\begin{proof}
Put $M=C_{p^8}$ and
$V=\{(a_i)\in M^q:\sum_i a_i=0\}$. Cyclic permutation of coordinates
has order $q$ and fixes only zero on $V$. As $V\cong M^{q-1}$, the
group $H\times V$ is a direct product of the IYB group $H\times M$
and the abelian group $M^{q-2}$, and hence is IYB by taking the direct product of their bijective cocycles.
The fixed subgroup for the permutation acting trivially on $H$ is
$H$. Any brace addition preserved by this automorphism would restrict to a brace addition on $H$, contradicting Theorem~\ref{thm:main}.
\end{proof}


\begin{thebibliography}{9}
\bibitem{Bac} D.~Bachiller, Counterexample to a conjecture about braces,
\emph{J. Algebra} \textbf{453} (2016), 160--176.
\bibitem{CJR} F.~Ced\'o, E.~Jespers and \'A.~del R\'io,
Involutive Yang--Baxter groups, \emph{Trans. Amer. Math. Soc.}
\textbf{362} (2010), 2541--2558.
\bibitem{CO} F.~Ced\'o and J.~Okni\'nski, New classes of IYB groups,
arXiv:2502.14805 (2025). \url{https://arxiv.org/abs/2502.14805}.
\bibitem{Eis} F.~Eisele, On the IYB-property in some solvable groups,
\emph{Arch. Math. (Basel)} \textbf{101} (2013), 309--318.
\bibitem{KN} E.~I.~Khukhro and V.~D.~Mazurov (eds.),
\emph{The Kourovka Notebook}, 21st issue (2026), Problem~21.7,
proposed by A.~Ballester-Bolinches. \url{https://alglog.org/21tkt.pdf}.
\bibitem{MBJ} H.~Meng, A.~Ballester-Bolinches and P.~Jim\'enez-Seral,
On a class of involutive Yang--Baxter groups, \emph{J. Algebra}
\textbf{704} (2026), 1--16.
\bibitem{Tr} S.~Trappeniers, A Lazard correspondence for post-Lie rings
and skew braces. \url{https://arxiv.org/abs/2406.02475}.
\bibitem{TrThesis} S.~Trappeniers, \emph{Skew braces and their connection
to the Yang--Baxter equation, Hopf--Galois structures and post-Lie rings},
Ph.D. thesis, Vrije Universiteit Brussel, Corollary~8.7.5.
\url{https://leandrovendramin.org/files/trappeniers.pdf}.
\bibitem{Ven} L.~Vendramin, Problems on skew left braces,
\emph{Adv. Group Theory Appl.} \textbf{7} (2019), 15--37.
\end{thebibliography}
\end{document}